# OPTIMAL SCALING OF MALA FOR NONLINEAR REGRESSION[1]

By Laird Arnault Breyer, Mauro Piccioni and Sergio Scarlatti

*University of Lancaster, University of Rome La Sapienza and University G. D'Annunzio Chieti*

We address the problem of simulating efficiently from the posterior distribution over the parameters of a particular class of nonlinear regression models using a Langevin–Metropolis sampler. It is shown that as the number $N$ of parameters increases, the proposal variance must scale as $N^{-1/3}$ in order to converge to a diffusion. This generalizes previous results of Roberts and Rosenthal [*J. R. Stat. Soc. Ser. B Stat. Methodol.* **60** (1998) 255–268] for the i.i.d. case, showing the robustness of their analysis.

**1. Introduction.** The motivation for the study of the kind of models analyzed in the present paper is the following. We consider a sequence of nonlinear regression models (indexed by $N$) relating a scalar response variable $y$ with a vector of covariates $z$

$$(1) \qquad y = \frac{1}{N} \sum_{i=1}^{N} h(z; x_i) + \frac{\varepsilon}{\sqrt{N}},$$

where $h(\cdot; x)$ is some function depending on a $d$-dimensional vector of parameters $x$ (weights) and $\varepsilon$ has a standard Gaussian distribution. If we take $n$ independent measurements $\mathbf{Y} = (Y_1, \ldots, Y_n)$ on the response variable, corresponding to the values $(z_1, \ldots, z_n)$ for the covariates, and define the vector $\mathbf{H}$ with components $H_k(x) = h(z_k; x)$, $k = 1, \ldots, n$, we get the measurement equation

$$(2) \qquad \mathbf{Y} = \frac{1}{N} \sum_{i=1}^{N} \mathbf{H}(x_i) + \frac{\varepsilon}{\sqrt{N}},$$

Received October 2000; revised May 2003.
[1]Supported by a TMR-grant "Spatial and Computational Statistics."

where $\boldsymbol{\varepsilon} = (\varepsilon_1, \ldots, \varepsilon_n)$ is a vector of i.i.d. standard Gaussians.

Following the Bayesian approach we take the vector of weights $(X_1, \ldots, X_N)$ to be random with i.i.d. $\mu$ distributed components. Then the measurement equation induces the following posterior distribution (i.e., conditional on $\mathbf{Y} = \mathbf{y}$) on the weights

$$\pi_N(dx) = C_N^{-1} \exp\left(\sum_{i=1}^N \langle \mathbf{y}, \mathbf{H}(x_i)\rangle - \frac{1}{2N} \sum_{i,j=1}^N \langle \mathbf{H}(x_i), \mathbf{H}(x_j)\rangle \right) \bigotimes_{i=1}^N \mu(dx_i),$$
(3)

where $\langle \cdot, \cdot \rangle$ stands for the usual scalar product in $\mathbb{R}^n$.

These kind of distributions are known in the statistical mechanics setting as "mean field" models [12]. The study of such distributions with a general nonlinear $\mathbf{H}$ is made complicated by the interaction term which destroys the a priori independence among the weights. In Appendix A we recall that propagation of chaos holds for the sequence of distributions (3) as $N \to \infty$ (Proposition 3, see also [1, 9]), which means that in the limit any finite collection of variables behaves as if the individual components had been drawn independently from a single probability measure $\pi$. This is characterized by

$$\log(d\pi/d\mu)(x) \propto \left\langle \mathbf{y} - \int \mathbf{H}\, d\pi, \mathbf{H}(x) \right\rangle.$$

Moreover, we prove a moderate deviations result (Proposition 5) which will be useful for the sequel.

In the rest of the paper we shall analyze the behavior of the Metropolis-adjusted Langevin algorithm (MaLa) [16] for distributions of the type (3). In order to simplify our analysis we shall consider the simplest case in which $n = 1$ and the weights are one-dimensional. Moreover, we shall assume that $\mu$ has an everywhere positive density w.r.t. the Lebesgue measure so the measure (3) has in this case the following $N$-dimensional posterior density

$$\pi_N(x) \propto \exp\left(\sum_{i=1}^N U(x_i) - \frac{1}{2N} \sum_{i,j=1}^N H(x_i) H(x_j)\right),$$
(4)

where

$$U(x) = yH(x) + \log \frac{d\mu}{dx}(x)$$

and the limiting probability measure $\pi$ on the real line has a positive density as well (called again $\pi$ to keep the notation simpler) with the property

$$\log \pi(x) \propto U(x) - H(x) \int H\, d\pi =: \psi(x).$$

In the following $X$ will always denote a random variable with density $\pi$ and expected values of measurable functions $f(X)$ will be written as $\pi(f(X))$.

The MaLa for the above density is a Markovian algorithm implemented in the following way. In order to compute $X_{j+1}^{(N)}$ given $X_j^{(N)}$, first generate

$$Y_j^{(N)} = X_j^{(N)} + \sigma W + \frac{\sigma^2}{2} \nabla \log \pi_N(X_j^{(N)}), \tag{5}$$

where $W$ is a standard Gaussian on $\mathbb{R}^N$ independent of $X_j^{(N)}$. The law of $Y_j^{(N)}$ given $X_j^{(N)} = x$, thus, has the density

$$\begin{aligned} q_N(x,y) \\ &\propto \exp\left(-\frac{1}{2\sigma^2}\left\|y - x - \frac{\sigma^2}{2}\nabla \log \pi_N(x)\right\|^2\right) \\ &= \exp\left(-\frac{1}{2\sigma^2}\sum_{i=1}^N \left(y_i - x_i - \frac{\sigma^2}{2}U'(x_i) - \frac{1}{N}\sum_{j=1}^N H'(x_i)\cdot H(x_j)\right)^2\right). \end{aligned} \tag{6}$$

The proposal $Y_j^{(N)}$ is accepted or rejected according to the following rule:

$$X_{j+1}^{(N)} = \begin{cases} Y_j^{(N)}, & \text{if } \xi_{j+1} < \dfrac{\pi_N(Y_j^{(N)})q_N(Y_j^{(N)},X_j^{(N)})}{\pi_N(X_j^{(N)})q_N(X_j^{(N)},Y_j^{(N)})}, \\ X_j^{(N)}, & \text{otherwise,} \end{cases} \tag{7}$$

where $\xi_j$ are i.i.d. $U[0,1]$.

In order to make the algorithm efficient the parameter $\sigma$ has to scale with $N$. A thorough discussion of this problem is reported in the recent survey [15], to which the reader is referred for more details. In the i.i.d. case ($H=0$), the optimal solution for the MaLa has been given by Roberts and Rosenthal [14]. Our main result is a generalization of theirs for sequences of densities of the type (4): if $\sigma$ is taken proportional to a suitable inverse power of the number of variables then the rescaled path of the algorithm converges weakly to a product of one-dimensional diffusions with the same stationary density $\pi(x)$. The choice of the proportionality factor only changes the (constant) speed at which the paths of the diffusions are travelled.

THEOREM 1 (Weak convergence of the MaLa). *Assume:*

(HP) *The functions $H$ and $U$ have bounded derivatives of all orders; moreover, $H$ itself is bounded, whereas $\lim_{|x|\to\infty} U(x) = -\infty$.*

Let $X_j^{(N)} = (X_j^{(N),1},\ldots,X_j^{(N),N})$ *be the MaLa defined by* (7), *with* $X_0^{(N)} \sim \pi_N$ *and* $\sigma^2 = \ell^2/N^{1/3}$. *The following weak convergence result holds in the*

*space $D[0,T]$,*

$$
\{(X^{(N),1}_{[tN^{1/3}]}, \ldots, X^{(N),k}_{[tN^{1/3}]}) : t \in [0,T]\}
$$
(8)
$$
\implies \{(Z^1_{v(\ell)t}, \ldots, Z^k_{v(\ell)t}) : t \in [0,T]\},
$$

*for any integer $k$, where $\{Z^i_t : i = 1, 2, \ldots\}$ are independent copies of the process $Z_t$ which is the unique solution to the SDE*

(9) $$dZ_t = \tfrac{1}{2}(\log \pi)'(Z_t)\,dt + dB_t, \qquad Z_0 \sim \pi,$$

*with $v = v(\ell) := 2\ell^2 \Phi(-\ell^3 \tau/2)$, $\tau$ being a constant depending on $\pi$ (explicitly given in Lemma 7 in Appendix B). Moreover, the acceptance probability converges as $N \to \infty$,*

$$
\lim_{N \to \infty} P(X^{(N)}_{j+1} = Y^{(N)}_j) = 2\Phi(-\ell^3 \tau/2) =: a(\ell).
$$

An implication of this result is that as $N \to \infty$, for any $T > 0$

(10) $$\frac{1}{TN^{1/3}} \sum_{j=1}^{TN^{1/3}} g(X^{(N)}_j) \to \frac{1}{v(\ell)T} \int_0^{v(\ell)T} g(Z_s)\,ds$$

weakly, if $g$ is bounded and continuous and depends only on $k$ components. Now, by the propagation of chaos, when $N$ is sufficiently large, the asymptotic bias

$$
\int g(x_1, \ldots, x_k) \pi_N(x_1, \ldots, x_N)\,dx_1 \cdots dx_N
$$
$$
- \int g(x_1, \ldots, x_k) \pi(x_1)\,dx_1 \cdots \pi(x_k)\,dx_k
$$

is small. On the other hand, by ergodicity of (9), when $T$ is large enough the right-hand side of (10) will be close to $\int g(x_1, \ldots, x_k) \pi(x_1) \cdots \pi(x_k)\,dx_1 \cdots dx_k$ with arbitrarily high probability [see, e.g., [17], Theorem (53.1)]. Hence, (10) may be loosely interpreted as stating that the Monte Carlo estimate

(11) $$\frac{1}{I} \sum_{j=1}^{I} g(X^{(N),1}_j, \ldots, X^{(N),k}_j)$$

of $\int g(x_1, \ldots, x_k) \pi_N(x_1, \ldots, x_N)\,dx_1 \cdots dx_N$ requires a number of iterations $I$ proportional to $N^{1/3}$. How large $T$ must be depends on the mixing properties of the diffusion $Z$, but it is, however, clear that for any fixed value of $T$ it is convenient to have $v(\ell)$ as large as possible in order to enlarge as much as possible the integration window. We can give an analytic expression for the maximizer $\hat{\ell}$ of $v(\ell)$, but this is, in practice, useless since it cannot be computed easily (except by Monte Carlo methods, which defeats somewhat the purpose). Luckily, the functions $v(\ell)$ and $a(\ell)$ have the same form as in [14], even if the constant $\tau$ is different in general. Hence, we can exploit the fact that $a$ is a bijective function of $\ell$ in order to maximize easily $v$ as a function of $a$. Indeed, $v(a) \propto a\{\Phi^{-1}(a/2)\}^{2/3}$, up to a constant factor depending on $\tau$. Since this function has a unique maximum in $a \approx 0.574$, in practice it suffices to monitor the acceptance rate $\frac{1}{k}\sum_{j=1}^{k} \mathbb{1}\{X_{j+1}^{(N)} \neq X_j^{(N)}\}$ of the MaLa and tune $\ell$ until $a(\ell)$ equals $0.574$.

As in the i.i.d. case, it is worth noticing the superiority of the MaLa over the random walk Metropolis (RWM) algorithm. In the RWM algorithm the proposal vector $Y^{(N)}$ has zero mean and, in order to obtain convergence to a diffusion $N^{1/3}$ has to be replaced by $N$, both in the scaling for the variance and for the time. The original result in [13] has been extended in [2] to Gibbs fields with no phase transition, and it could be proved for mean field models like (4) as well. As a consequence, (10) essentially holds with $N^{1/3}$ replaced by $N$, which implies that the required number of steps has the order $N$ rather than $N^{1/3}$. The only difference is that the function $v(\ell)$ has to be replaced by some other function, which this time is maximized when the acceptance rate is roughly equal to $0.234$.

A final comment concerns the assumption made in Theorem 1 that the initial value $X_0^{(N)}$ is already distributed according to the target density $\pi_N$, which is clearly unrealistic. This means that, in practice, the partial sums in (10) do not start from 1, but typically from some large value $t_0$, which ensures that the effect of the initial value $X_0^{(N)}$ can be neglected. A deeper study of the scaling behavior of the MaLa and the RWM when started in the tails of the target density $\pi_N$ has been initiated in [3].

**2. A quantitative central limit theorem for the log-acceptance ratio.** A fundamental step towards the proof of Theorem 1 is to establish a quantitative central limit theorem (CLT) for the log-acceptance ratio

$$(12) \qquad G_{\sigma,N}(x,W) = \log \frac{\pi_N(Y_\sigma(x,W))q_N(Y_\sigma(x,W),x)}{\pi_N(x)q_N(x,Y_\sigma(x,W))},$$

where $x = (x_1, \ldots, x_N)$ is fixed, $W = (W_1, \ldots, W_N)$ is a random vector having i.i.d. $N(0,1)$ components defined on some probability space $(\Omega, \mathcal{F}, \mathbb{P})$ and $Y_\sigma(x,W)$ is the proposal vector given by

$$(13) \quad Y_{\sigma,i}(x,W) = Y_i = x_i + \sigma W_i + \frac{\sigma^2}{2}\left(U'(x_i) - H'(x_i)\frac{1}{N}\sum_{j=1}^{N} H(x_j)\right),$$

for $i = 1, \ldots, N$, with $\sigma = \sigma_N = \frac{\ell}{N^{1/6}}$, for some $\ell > 0$.

PROPOSITION 2 (CLT for the acceptance ratio). *There exist measurable sets $F_N \subset \mathbb{R}^N$, with $\pi_N(F_N^c) = o(N^{-t})$ for any $t > 0$, such that*

$$(14) \quad \lim_{N \to \infty} N^\beta \sup_{x \in F_N} \sup_{u \in R} \left| \mathbb{P}\left( \frac{G_{\sigma_N, N}(x, W)}{\ell^3 \tau} + \frac{\ell^3 \tau}{2} \leq u \right) - \Phi(u) \right| = 0$$

*for any $\beta > 0$ sufficiently small, where $\tau$ is some positive constant.*

Before starting the proof we set up a convenient notation. First, we shall denote by $E_N$ empirical averages w.r.t. the vector $(x, W, Y)$, that is,

$$(15) \quad E_N f(x, W, Y) = \frac{1}{N} \sum_{i=1}^{N} f(x_i, W_i, Y_i).$$

In order to shorten the notation even further the function $f$ is allowed to contain empirical averages as arguments as well, in which case they have to be considered as constants. In particular, for

$$(16) \quad \psi_N(t; x) = U(t) - H(t) E_N H(x),$$

we define

$$(E_N \psi_N)(x) = \frac{1}{N} \sum_{i=1}^{N} \psi_N(x_i; x) = E_N U(x) - (E_N H(x))^2,$$

and we apply the same convention to empirical averages of derivatives

$$\psi_N^{(k)}(t; x) = U^{(k)}(t) - H^{(k)}(t) E_N H(x)$$

and to their products. Finally, we use the shortened notation

$$(17) \quad E_N g(x) W^l = \frac{1}{N} \sum_{i=1}^{N} g(x_i) W_i^l$$

and

$$(18) \quad E_N h(Y) W^l = \frac{1}{N} \sum_{i=1}^{N} h(Y_i) W_i^l.$$

Moreover, we will always use the same letter $C$ for several constants appearing in the estimates.

PROOF OF PROPOSITION 2. By direct computation the first two derivatives of $G_{\sigma, N}(x, W)$ w.r.t. $\sigma$ vanish at $\sigma = 0$. Consequently, we have the Taylor expansion

$$(19) \quad G_{\sigma, N}(x, W) = \sum_{k=3}^{6} \sigma^k g_{k, N}(x, W) + \frac{1}{6!} \int_0^\sigma (\sigma - u)^6 \frac{d^7}{du^7} G_{u, N}(x, W) \, du,$$

where $g_{k,N}(x,W) = \frac{1}{k!}\frac{d^k}{du^k}G_{u,N}(x,W)(0)$ for $k = 3,\ldots,6$. For completeness the explicit form of these functions is given in Lemma 6 in Appendix B. Setting $\sigma = \ell/N^{1/6}$ and standardizing as in (14), we have

$$\frac{G_{\sigma_N,N}}{\ell^3\tau} + \frac{\ell^3\tau}{2} = \frac{1}{N^{1/2}\tau}g_{3,N}(x,W) + \frac{\ell}{N^{2/3}\tau}g_{4,N}(x,W)$$
$$+ \frac{\ell^2}{N^{5/6}\tau}g_{5,N}(x,W) + \frac{\ell^3}{\tau}\left(g_{6,N}(x,W) + \frac{\tau^2}{2}\right)$$
$$+ \frac{1}{6!\tau\ell^3}\int_0^{\ell N^{-1/6}}(\ell N^{-1/6} - u)^6\frac{d^7}{du^7}G_{u,N}(x,W)\,du$$
$$=: A_N + B_N + C_N + D_N + I_N.$$

By using a a standard lemma on distribution functions ([11], Lemma 1.9, page 20) we obtain the following estimate:

$$\sup_{u\in R}\left|\mathbb{P}\left(\frac{G_{\sigma_N,N}}{\ell^3\tau} + \frac{\ell^3\tau}{2} \le u\right) - \Phi(u)\right|$$

(20)
$$\le \sup_{u\in R}|\mathbb{P}(A_N \le u) - \Phi(u)| + \mathbb{P}(|B_N| \ge \varepsilon_N) + \mathbb{P}(|C_N| \ge \varepsilon_N)$$
$$+ \mathbb{P}(|D_N| \ge \varepsilon_N) + \mathbb{P}(|I_N| \ge \varepsilon_N) + \frac{4\varepsilon_N}{\sqrt{2\pi}},$$

where $(\varepsilon_N)$ is an arbitrary sequence of positive numbers to be chosen in the sequel.

In Appendix B various lemmas are proven in order to estimate separately each term appearing on the right-hand side of (20). By Lemma 7, for any $N$ and $\varepsilon_N > 0$,

(21)
$$\sup_{u\in R}|\mathbb{P}(A_N \le u) - \Phi(u)|$$
$$\le C\left(\frac{1}{\sqrt{N}} + \frac{1}{\varepsilon_N^2 N}\right) + h_\tau(F_3(E_N\mathbf{r}_3(x))) + \frac{\varepsilon_N}{\sqrt{2\pi}},$$

where $F_3$ is polynomial, $\mathbf{r}_3$ is a vector of bounded measurable functions and $h_\tau$ is a locally Lipschitz function vanishing at

$$\tau^2 = F_3(\pi(\mathbf{r}_3(X))).$$

Denote by $C_3$ the inverse of the local Lipschitz constant of $h$ at $\tau^2$. Therefore, for

$$x \in F_{N,3}(\varepsilon_N) = \{x : |E_N\mathbf{r}_3(x) - \pi(\mathbf{r}_3(X))| \le C_3\varepsilon_N\},$$

it holds

(22) $$\sup_{u\in R}|\mathbb{P}(A_N \le u) - \Phi(u)| \le C\left(\frac{1}{\sqrt{N}} + \frac{1}{\varepsilon_N^2 N} + \varepsilon_N\right),$$

provided $\varepsilon_N$ goes to zero. By Lemma 11, for any $N$ and $\varepsilon_N > 0$,

$$\mathbb{P}(|B_N| \geq \varepsilon_N) \leq \frac{C}{N^{1/3}\varepsilon_N^2}, \tag{23}$$

$$\mathbb{P}(|C_N| \geq \varepsilon_N) \leq \frac{C}{N^{2/3}\varepsilon_N^2}, \tag{24}$$

$$\mathbb{P}(|D_N| \geq \varepsilon_N) \leq \frac{C}{N\varepsilon_N^2}, \tag{25}$$

for $x \in \bigcap_{k=4}^{6} F_{N,k}(\varepsilon_N)$, where

$$F_{N,k}(\varepsilon_N) = \{x : |E_N \mathbf{r}_k(x) - \pi(\mathbf{r}_k(X))| \leq C_k \varepsilon_N N^{k/6-1}\},$$

$\mathbf{r}_k$ being a vector of functions for $k = 4, 5, 6$, and $C_k$, $k = 4, 5, 6$, are suitably small constants. Furthermore, by Lemma 12,

$$\mathbb{P}(|I_N| \geq \varepsilon_N) \leq \frac{C}{N^{1/6}\varepsilon_N}. \tag{26}$$

Finally, set $F_N = \bigcap_{k=3}^{6} F_{N,k}(\varepsilon_N)$, and choose $\varepsilon_N = N^{-1/9}$. In order to estimate $\pi_N(F_{N,k}^c(N^{-1/9}))$ we need to control deviations of empirical averages from expected values under $\pi$ of the order $N^{-\alpha_k}$, where $\alpha_3 = \alpha_6 = 1/9$, $\alpha_5 = 5/18$ and $\alpha_4 = 4/9$. Since the latter is the largest, it is enough to apply Proposition 5 in Appendix A with $\lambda_N = N^{1/18}$, in which case $N^{-1/2}\lambda_N = N^{-4/9}$. By consequence,

$$\pi_N(F_N^c) \leq \sum_{k=3}^{6} \pi_N(F_{N,k}(N^{-1/9})) \leq \exp(-cN^{1/9} + o(N^{1/9})),$$

which is $o(N^{-t})$ for any $t > 0$ as claimed.

Using the bounds (20), (22)–(26) we get that

$$\sup_{u \in R} |\mathbb{P}(G_{\sigma,N}(x,W) \leq u) - \Phi_{-\ell^6\tau^2/2,\ell^6\tau^2}(u)| = O(N^{-1/9}). \quad \square$$

**3. Proof of Theorem 1.** Let $f$ be any smooth function with compact support from $\mathbb{R}^N$ to $\mathbb{R}$. Define on $f$ the discrete generator,

$$\begin{aligned} A_{\sigma,N}f(x) &= \mathbb{E}[f(X_{t+1}^{(N)}) - f(x)|X_t^{(N)} = x] \\ &= \mathbb{E}[(f(Y_\sigma) - f(x))1 \wedge e^{G_{\sigma,N}(x,W)}], \end{aligned} \tag{27}$$

and the infinitesimal generator of the process $(Z_{v(\ell)t})$,

$$Af(x) = \frac{v(\ell)}{2} \sum_{p=1}^{N} \left[ f_{x_p x_p}(x) + \left( U'(x_p) - H'(x_p) \int H \, d\pi \right) f_{x_p}(x) \right]. \tag{28}$$

By [7], Corollary 8.9, page 233, the weak convergence (8) holds, provided we exhibit measurable sets $\widetilde{F}_N \subset \mathbb{R}^N$ such that

$$\lim_{N \to \infty} \mathbb{P}(X^{(N)}_{[N^{1/3}t]} \in \widetilde{F}_N \text{ for all } t \leq T) = 1 \tag{29}$$

and

$$\lim_{N \to \infty} \sup_{x \in \widetilde{F}_N} |N^{1/3} A_{\ell N^{-1/6}, N} f(x) - Af(x)| = 0 \tag{30}$$

for any smooth $f(x) = f(x_1, \ldots, x_k)$ with compact support. Notice that since $X_0^{(N)} \sim \pi_N$ and $\pi_N$ is stationary,

$$\mathbb{P}(X^{(N)}_{[N^{1/3}t]} \notin \widetilde{F}_N \text{ for some } t \leq T) \leq [N^{1/3}T]\pi_N(x : x \notin \widetilde{F}_N).$$

Thus, in order to ensure (29) it is enough to check that $\pi_N(\widetilde{F}_N^c) = o(N^{-1/3})$. By [18], Proposition 2.2, page 177, it is enough to prove (30), for $k = 2$, in order to get the convergence (8) for any integer $k$.

For a fixed $x \in \mathbb{R}^N$ we expand $A_{\sigma,N} f(x)$ in powers of $\sigma$, which is obtained by recalling that $Y_{\sigma,1}$ is defined in (13),

$$\begin{aligned}
A_{\sigma,N} f(x) &= \mathbb{E}[(f(Y_{\sigma,1}, Y_{\sigma,2}) - f(x_1, x_2)) 1 \wedge e^{G_{\sigma,N}(x,W)}] \\
&= \mathbb{E}\bigg\{\bigg[\sum_{i=1}^{2}\bigg(\sigma W_i f_{x_i} + \frac{\sigma^2}{2} W_i^2 f_{x_i x_i} \\
&\qquad + \frac{\sigma^2}{2} f_{x_i}(U'(x_i) - H'(x_i)E_N H(x))\bigg) + \sigma^2 W_1 W_2 f_{x_1 x_2}\bigg] \\
&\qquad \times \mathbb{E}(1 \wedge e^{G_{\sigma,N}(x,W)} | W_1, W_2)\bigg\} \\
&\quad + \sigma^3 r_N(\sigma, x),
\end{aligned} \tag{31}$$

where partial derivatives of $f$ are always evaluated at $(x_1, x_2)$ if not specified otherwise, and

$$\begin{aligned}
r_N(\sigma, x) = \frac{\sigma}{3!} \mathbb{E}\bigg\{\bigg(\sum_{i=1}^{2}\bigg[&f_{x_i, x_i, x_i}(Y_{\widetilde{\sigma},1}, Y_{\widetilde{\sigma},2}) \\
&\times (W_i + \widetilde{\sigma}(U'(x_i) - H'(x_i)E_N H(x)))^3 \\
&+ 3 f_{x_i, x_i}(Y_{\widetilde{\sigma},1}, Y_{\widetilde{\sigma},2}) \cdot (U'(x_i) - H'(x_i)E_N H(x)) \\
&\times (W_i + \widetilde{\sigma}(U'(x_i) - H'(x_i)E_N H(x)))\bigg] \\
+ 3 \sum_{i \neq j} \bigg[&f_{x_i, x_i, x_j}(W_i + \widetilde{\sigma}(U'(x_i) - H'(x_i)E_N H(x)))^2
\end{aligned}$$

$$\times (W_j + \tilde{\sigma}(U'(x_j) - H'(x_j)E_N H(x)))$$
$$+ f_{x_i,x_j}(W_i + \tilde{\sigma}(U'(x_i) - H'(x_i)E_N H(x)))$$
$$\times (U'(x_i) - H'(x_i)E_N H(x))\Big]\Big)$$
$$\times 1 \wedge e^{G_{\sigma,N}(x,W)}\Big\},$$

where $0 \leq \tilde{\sigma} \leq \sigma$. By assumption (HP), plugging in $\sigma = \sigma_N = \ell N^{-1/6}$, the remainder $r_N(\sigma, x)$ is uniformly bounded in $N$ and $x$.

Next, observe that if $\Gamma(u)$ is an absolutely continuous function of the real variable $u$, then

$$1 \wedge e^{\Gamma(1)} = 1 \wedge e^{\Gamma(0)} + \int_0^1 \mathbb{1}_{\{\Gamma(u)<0\}} \Gamma'(u) e^{\Gamma(u)} \, du.$$

Now we apply this formula to the function $\tilde{\Gamma}(u) = G_{\sigma,N}(x, uW_1, uW_2, W^{(c)})$, where $W^{(c)} = (W_3, \ldots, W_N)$, and take conditional expectations w.r.t. $(W_1, W_2)$:

$$\mathbb{E}(1 \wedge e^{G_{\sigma,N}(x,W)} | W_1, W_2)$$
$$(32) \quad = \mathbb{E}(1 \wedge e^{G_{\sigma,N}(x,W)} | W_1 = 0, W_2 = 0)$$
$$+ \sum_{i=1}^{2} W_i \int_0^1 \mathbb{E}(\mathbb{1}_{\{\tilde{\Gamma}(u)<0\}} G^{(i)}_{\sigma,N}(x, uW_1, uW_2, W^{(c)}) e^{\tilde{\Gamma}(u)} | W_1, W_2) \, du,$$

where $G^{(i)}_{\sigma,N}$ denotes the partial derivative of $G_{\sigma,N}(x, w_1, w_2, w^{(c)})$ w.r.t. the variable $w_i$.

We now substitute (32) into the expression (31) so we obtain the following expression:

$$\frac{1}{\sigma^2} A_{\sigma,N} f(x) = \frac{1}{2} \sum_{i=1}^{2} \Big[ f_{x_i,x_i} \mathbb{E}(1 \wedge e^{G_{\sigma,N}(x,W)} | W_1 = 0, W_2 = 0)$$
$$(33) \qquad + f_{x_i}(U'(x_i) - H'(x_i)E_N H(x)) \mathbb{E}(1 \wedge e^{G_{\sigma,N}(x,W)}) \Big]$$
$$+ R_N(\sigma, x),$$

where

$$R_N(\sigma, x)$$
$$= \sigma^{-1} \sum_{i=1}^{2} f_{x_i} \mathbb{E}\Big[ W_i^2 \int_0^1 \mathbb{1}_{\{\tilde{\Gamma}(u)<0\}} G^{(i)}_{\sigma,N}(x, uW_1, uW_2, W^{(c)}) e^{\tilde{\Gamma}(u)} \, du \Big]$$

$$+ \frac{1}{2}\sum_{i=1}^{2} f_{x_i,x_i} \mathbb{E}\left[W_i^2 \int_0^1 \mathbb{1}_{\{\tilde{\Gamma}(u)<0\}} G_{\sigma,N}^{(i)}(x, uW_1, uW_2, W^{(c)}) e^{\tilde{\Gamma}(u)}\, du\right]$$

$$+ f_{x_1,x_2}\left\{\mathbb{E}\left[W_1^2 W_2 \int_0^1 \mathbb{1}_{\{\tilde{\Gamma}(u)<0\}} G_{\sigma,N}^{(1)}(x, uW_1, uW_2, W^{(c)}) e^{\tilde{\Gamma}(u)}\, du\right]\right.$$

$$\left.+ \mathbb{E}\left[W_1 W_2^2 \int_0^1 \mathbb{1}_{\{\tilde{\Gamma}(u)<0\}} G_{\sigma,N}^{(2)}(x, uW_1, uW_2, W^{(c)}) e^{\tilde{\Gamma}(u)}\, du\right]\right\}$$

$$+ \sigma r_N(\sigma, x).$$

Now we concentrate on the $\sigma^{-1}$ term in the above expression, since the others are more easily controlled with similar arguments. First, bound $|f_{x_i}|$ with a constant, then we are left to bound for $i = 1, 2$,

$$\text{(34)} \quad \begin{aligned} &\frac{1}{\sigma}\mathbb{E}\left[W_i^2 \int_0^1 \mathbb{1}_{\{\tilde{\Gamma}(u)<0\}} G_{\sigma,N}^{(i)}(x, uW_1, uW_2, W^{(c)}) e^{\tilde{\Gamma}(u)}\, du\right] \\ &\leq \frac{1}{\sigma}\mathbb{E}\left(W_i^2 \sup_{0\leq u\leq 1} |G_{\sigma,N}^{(i)}(x, uW_1, uW_2, W^{(c)})|\right). \end{aligned}$$

Let us write explicitly

$$G_{\sigma,N}^{(i)}(x, W)$$
$$= \frac{\sigma}{2}(U'(Y_i) - U'(x_i) - H'(Y_i)E_N H(Y) + H'(x_i)E_N H(x))$$
$$- \frac{\sigma^2}{2}\left(W_i(U''(Y_i) - H''(Y_i)E_N H(Y)) - H'(Y_i)\frac{1}{N}\sum_{k=1}^{N} W_k H'(Y_k)\right)$$
$$+ \frac{\sigma^3}{8}(U'''(Y_i) - H'''(Y_i)E_N H(Y) - H'(Y_i)E_N H'(Y)),$$

where we have written $Y_i$ for $Y_{\sigma,i}$. Using (HP), we can rewrite the right-hand side of (34) as

$$\frac{1}{\sigma}\mathbb{E}\left(W_i^2 \sup_{0\leq u\leq 1} |G_{\sigma,N}^{(i)}(x, uW_1, uW_2, W^{(c)})|\right)$$
$$= \frac{1}{2}\mathbb{E}\left(W_i^2 \sup_{0\leq u\leq 1} |U'(Y_i(u)) - U'(x_i)\right.$$
$$\left. - H'(Y_i(u))E_N H(Y(u)) + H'(x_i)E_N H(x)|\right) + o(1),$$

where $Y_k(u) = Y_k + \sigma u \sum_{i=1}^{2} \delta_{ki} W_i$, $k = 1, \ldots, N$. We have now

$$\mathbb{E}\left(W_i^2 \mathbb{E}\left(\sup_{0\leq u\leq 1} |U'(Y_i(u)) - U'(x_i)\right.\right.$$

$$- H'(Y_i(u))E_N H(Y(u)) + H'(x_i)E_N H(x)|\Big)\Big)$$

(35)
$$\leq \mathbb{E}\Big(W_i^2 \sup_{0\leq u\leq 1} |U'(Y_i(u)) - U'(x_i)|\Big)$$

$$+ \mathbb{E}\Big(W_i^2 \sup_{0\leq u\leq 1} |H'(x_i) - H'(Y_i(u))|E_N H(x)\Big)$$

$$+ \mathbb{E}\Big(W_i^2 \sup_{0\leq u\leq 1} (|H'(Y_i(u))|E_N|H(Y(u)) - H(x)|)\Big).$$

Observe that, when $T$ is either $U'$, $H'$ or $H$ and $i = 1, 2$, we can write, using the fundamental theorem of calculus,

$$T(Y_i(u)) - T(x_i)$$
$$= T\Big[x_i + \frac{\sigma^2}{2}(U'(x_i) - H'(x_i)E_N H(x))\Big] - T(x_i) + T(Y_i(u))$$
$$\quad - T\Big[x_i + \frac{\sigma^2}{2}(U'(x_i) - H'(x_i)E_N H(x))\Big]$$
$$= \frac{\sigma^2}{2}(U'(x_i) - H'(x_i)E_N H(x))$$
$$\quad \times \int_0^1 T'\Big(x_i + \frac{v\sigma^2}{2}(U'(x_i) - H'(x_i)E_N H(x))\Big)\, dv$$
$$\quad + \sigma W_i \int_0^u T'\Big(x_i + \frac{\sigma^2}{2}(U'(x_i) - H'(x_i)E_N H(x)) + s\sigma W_i\Big)\, ds.$$

By bounding the derivative of $T$ and substituting $\sigma_N = \ell N^{-1/6}$, we have

$$\sup_{0\leq u\leq 1} |T(Y_i(u)) - T(x_i)| \leq C N^{-1/6}(1 + |W_i|).$$

By substituting this bound into (35) and, subsequently in (34), the right-hand side is bounded by $O(N^{-1/6})$ uniformly over $x$. Similar arguments allow us to conclude that $R_N(\sigma, x) \to 0$ as $N \to \infty$, uniformly over $x$ as well.

Now let $\mathcal{N}$ be a Gaussian random variable with mean $-\ell^6\tau^2/2$ and variance $\ell^6\tau^2$. It is immediately seen that $\mathbb{E}(1 \wedge e^{\mathcal{N}}) = 2\Phi(-\ell^3\tau/2)$. By an integration by parts we have

$$|\mathbb{E}(1 \wedge e^{G_{\sigma,N}(x,W)}) - \mathbb{E}(1 \wedge e^{\mathcal{N}})|$$
$$\leq C \sup_{u\in R}|\mathbb{P}(G_{\sigma,N}(x,W) \leq u) - \Phi_{-\ell^6\tau^2/2,\ell^6\tau^2}(u)|,$$

which goes to zero uniformly for $x \in F_N$ by Lemma 7. Moreover,

$$|\mathbb{E}(1 \wedge e^{G_{\sigma,N}(x,W)}|W_1 = 0, W_2 = 0) - \mathbb{E}(1 \wedge e^{G_{\sigma,N}(x,W)})|$$

$$\leq \mathbb{E}|G_{\sigma,N}(x,0,0,W^{(c)}) - G_{\sigma,N}(x,W)|$$

$$\leq \sum_{i=1}^{2} \mathbb{E} \int_0^1 |W_i G_{\sigma,N}{}^{(i)}(x, uW_1, uW_2, W^{(c)})|\, du,$$

and by the same argument as before, the right-hand side goes to zero uniformly over $x$. Finally, we have

$$|N^{1/3} A_{\sigma,N} f(x) - Af(x)|$$
$$= \ell^2 |\sigma_N^{-2} A_{\sigma,N} f(x) - \ell^{-2} Af(x)|$$
$$\leq \tfrac{1}{2}\ell^{-2} \sum_{i=1}^{2} \Big[|f_{x_i,x_i}(x_1,x_2)||\mathbb{E}(1 \wedge e^{G_{\sigma,N}(x,W)}|W_1=0, W_2=0)$$
$$- \mathbb{E}(1 \wedge e^{\mathcal{N}})|$$
$$+ |f_{x_i}(x_1,x_2)||U'(x_i)||\mathbb{E}(1 \wedge e^{G_{\sigma,N}(x,W)}) - \mathbb{E}(1 \wedge e^{\mathcal{N}})|$$
$$+ |f_{x_i}(x_1,x_2)||H'(x_i)|$$
$$\times \Big(|E_N H(x)||\mathbb{E}(1 \wedge e^{G_{\sigma,N}(x,W)}) - \mathbb{E}(1 \wedge e^{\mathcal{N}})|$$
$$+ |E_N H(x) - \pi(H(X))||\mathbb{E}(1 \wedge e^{\mathcal{N}})|\Big)\Big]$$
$$+ |R_N(x,f)|.$$

Next define $\widetilde{F}_N = F_N \cap \{x : |E_N H(x) - \pi(H(X))| \leq N^{-1/9}\}$. By using Proposition 5 in Appendix A it is immediately verified that $\pi_N(\widetilde{F}_N^c) = o(N^{-t})$ for any $t > 0$. The proof is complete since the right-hand side of the last expression goes to zero uniformly on $\widetilde{F}_N$.

## APPENDIX A.

In this appendix we discuss the asymptotic behavior of sequences of distributions $\pi_N$ defined in (3) for a general measurable function $\mathbf{H}$. For ease of notation we drop from now on boldfaces used to indicate $n$-dimensional vectors. First let us introduce the exponential family of probability measures on $\mathbb{R}^d$ generated by $\mu$ and $H$, which is defined by

$$\mu_\theta(dx) = e^{\langle\theta, H(x)\rangle - K(\theta)} \mu(dx), \qquad \theta \in \Theta,$$

where $K(\theta) = \log \int e^{\langle\theta, H(x)\rangle} \mu(dx)$ is the cumulant generating function of $H$ under $\mu$. We assume that $K$ is finite only in an open set $\Theta$ of $\mathbb{R}^n$ and that no hyperplane of $\mathbb{R}^n$ contains $H(x)$ $\mu$-almost surely (in the case $n=1$ this is equivalent to assume that $H$ is nonconstant). Moreover, in the paper we assumed that $H$ is bounded so $K$ is defined on the whole space.

Consider now the strictly convex function $J(\theta) = \frac{1}{2}\|\theta - y\|^2 + K(\theta)$, where $K$ is extended to the complement of $\Theta$ by setting its value equal to $+\infty$. This function has a unique minimum $\theta_* = \theta_*(y)$ in $\mathbb{R}^n$ (as it is strictly convex and lower semicontinuous with compact level sets), that is the unique solution of the equation

$$\theta + \nabla K(\theta) = y, \tag{36}$$

which implies, by the properties of exponential families, that

$$\theta_* = y - \int H \, d\pi. \tag{37}$$

We can now state the following:

PROPOSITION 3 (Propagation of chaos). *Whenever $f : (\mathbb{R}^d)^\infty \to \mathbb{R}$ is a bounded measurable local function (i.e., it depends only on a finite number of components), then*

$$\lim_{N \to \infty} \int f \, d\pi_N = \int f \, d\pi^{\otimes \infty},$$

*where $\pi = \mu_{\theta_*}$.*

PROOF. We can easily bound the Kullback–Leibler divergence

$$D(\pi_N \| \pi^{\otimes N}) = \int \log(d\pi_N / d\pi^{\otimes N}) \, d\pi_N.$$

In fact, by using (37) and setting $\widetilde{H}(x) = H(x) - \int H \, d\pi$,

$$\log\left(\frac{d\pi_N}{d\pi^{\otimes N}}\right)$$

$$= \log C_N^{-1} + NK(\theta_*) + \sum_{i=1}^{N} \langle m - \theta_*, H(x_i)\rangle - \frac{1}{2N} \sum_{i,j=1}^{N} \langle H(x_i), H(x_j)\rangle$$

$$= \log C_N^{-1} + NK(\theta_*) + \frac{N}{2}\left\|\int H \, d\pi\right\|^2$$

$$- \frac{N}{2}\left(\left\|\int H \, d\pi\right\|^2 + \left\langle \sum_{i,j=1}^{N} \frac{H(x_i)}{N}, \frac{H(x_j)}{N}\right\rangle - 2\sum_{i=1}^{N}\left\langle \frac{H(x_i)}{N}, \int H \, d\pi\right\rangle\right)$$

$$= \log \widetilde{C}_N - \frac{1}{2}\left\|\frac{1}{\sqrt{N}}\sum_{i=1}^{N} \widetilde{H}(x_i)\right\|^2,$$

where

$$\log \widetilde{C}_N = \log C_N^{-1} + NK(\theta_*) + \frac{N}{2}\left\|\int H\,d\pi\right\|^2$$

$$= \log C_N^{-1} + NJ(\theta_*)$$

$$= -\log \int \exp\left(-\frac{1}{2}\left\|\frac{1}{\sqrt{N}}\sum_{i=1}^{N}\widetilde{H}(x_i)\right\|^2\right)\bigotimes_{i=1}^{N}\pi(dx_i),$$

and, therefore, by the CLT, the right-hand side of the above expression converges to

$$-\log E(\exp(-\tfrac{1}{2}|Z^2|)),$$

where $Z$ is a zero mean Gaussian vector. Hence, it is bounded uniformly in $N$ by some constant $M_0$. By consequence $D(\pi_N\|\pi^{\otimes N}) \leq M_0$. It follows that if we denote by $\pi_{N,k}$ the marginal of $\pi_N$ for the first $k$ components, then an inequality of Csiszar [5] equation (2.11), page 772, yields

$$D(\pi_{N,k}\|\pi^{\otimes k}) \leq \frac{1}{[N/k]}D(\pi_N\|\pi^{\otimes N}) \leq \frac{M_0}{[N/k]},$$

and now the stated convergence follows by [4], Lemma 3.1. □

In the forthcoming Proposition 5, we shall need the following technical lemma.

LEMMA 4. *For any symmetric nonnegative definite matrix $A$ of order $s$, the convex conjugate of $\theta \mapsto \frac{1}{2}\langle\theta, A\theta\rangle$ is given by*

$$M^*(z) = \begin{cases} \frac{1}{2}\langle z, A^- z\rangle, & \text{if } z \in \operatorname{Ran} A, \\ +\infty, & \text{otherwise,} \end{cases}$$

*where $A^-$ is the pseudo-inverse of $A$. As a consequence, the origin is the unique minimizer of $M^*$.*

PROOF. Let $A = U^t L U$, with $L$ a diagonal matrix with the diagonal elements equal to the eigenvalues $(\lambda_i)$ of $A$. Then $A^- = U^t L^- U$, where $L^-$ is the diagonal matrix with diagonal elements equal to the reciprocal of the eigenvalues (if positive) of $A$ and zero otherwise. By definition,

$$M^*(z) = \sup_{\theta}(\langle z, \theta\rangle - \tfrac{1}{2}\langle\theta, A\theta\rangle) = \sup_{w}\left(\sum_{i=1}^{s} v_i w_i - \tfrac{1}{2}\sum_{i=1}^{s}\lambda_i w_i^2\right),$$

where $v = Uz$ and $w = U\theta$. If there exists $i_0$ such that $\lambda_{i_0} = 0$ and $v_{i_0} \neq 0$ (which happens if and only if $z \notin \operatorname{Ran} A$), it is immediately seen that $M^*(z) = +\infty$. Otherwise, the function between round brackets has a maximum $w_i = \frac{v_i}{\lambda_i}$ for $i$ such that $\lambda_i > 0$, $w_i = 0$ otherwise. Finally, it is easily seen that

$$M^*(z) = \frac{1}{2} \sum_{i:\lambda_i>0} \frac{v_i^2}{\lambda_i} = \frac{1}{2}\langle z, A^- z \rangle$$

for $z \in \operatorname{Ran} A$. □

PROPOSITION 5 (Moderate deviations). *If the sequence $\{\lambda_N\}$ is such that $\lambda_N \to \infty$ but $\lambda_N^2/N \to 0$, then for any bounded measurable function $g : \mathbb{R}^d \to \mathbb{R}^m$,*

$$\pi_N\left(\left|\frac{1}{N}\sum_{i=1}^N g(x_i) - \int g\,d\pi\right| \geq \frac{\lambda_N}{\sqrt{N}}\right) \leq e^{-c\lambda_N^2 + o(\lambda_N^2)},$$

*where $c > 0$ is a constant and $\pi = \mu_{\theta_*}$.*

PROOF. Define $\widetilde{g}(x_i) = g(x_i) - \int g\,d\pi$, $\widetilde{H}(x_i) = H(x_i) - \int H\,d\pi$ and

$$(Z_N, Y_N) = (\lambda_N \sqrt{N})^{-1} \sum_{i=1}^N (\widetilde{g}(x_i), \widetilde{H}(x_i)).$$

Now it is easy to compute (see, e.g., [6])

$$\Lambda(\theta, \psi) = \lim_{N\to\infty} \frac{1}{\lambda_N^2} \log \int \exp \lambda_N^2 (\langle \theta, Z_N \rangle + \langle \psi, Y_N \rangle)\, d\pi^{\otimes N}$$

$$= \frac{1}{2}\langle(\theta,\psi), \Sigma(\theta,\psi)\rangle,$$

where $\Sigma$ is the covariance matrix of $(\widetilde{g}(x), \widetilde{H}(x))$ under $\pi$. By applying the Gärtner–Ellis theorem and Lemma 4, we prove that $(Z_N, Y_N)$ satisfies under $\pi$ an LDP with speed $\lambda_N^2$ and rate function

$$J(z,y) = \begin{cases} \frac{1}{2}\langle(z,y), \Sigma^-(z,y)\rangle, & \text{if } (z,y) \in \operatorname{Ran}\Sigma, \\ +\infty, & \text{otherwise.} \end{cases}$$

We want to prove the same result for the sequence $Z_N$ under $\pi_N$. Decompose $\Sigma$ into blocks as

$$\Sigma = \begin{pmatrix} \Sigma_{11} & \vdots & \Sigma_{12} \\ \cdots & & \cdots \\ \Sigma_{21} & \vdots & \Sigma_{22} \end{pmatrix} = \begin{pmatrix} \int \widetilde{g}\widetilde{g}^t\,d\pi & \vdots & \int \widetilde{g}\widetilde{H}^t\,d\pi \\ \cdots & & \cdots \\ \int \widetilde{H}\widetilde{g}^t\,d\pi & \vdots & \int \widetilde{H}\widetilde{H}^t\,d\pi \end{pmatrix}$$

and write

$$\widetilde{\Lambda}_N(\theta) = \log \widetilde{C}_N^{-1} \int \exp\{\lambda_N^2(\langle\theta, Z_N\rangle - \tfrac{1}{2}|Y_N|^2)\}\, d\pi^{\otimes N}$$

$$= \log \int \exp\{\lambda_N^2 \langle\theta, Z_N\rangle\}\, d\pi_N.$$

Next apply Varadhan's lemma ([6], Theorem 4.3.1, page 137) to the continuous function $\varphi(z,y) = \langle\theta, z\rangle - \tfrac{1}{2}\|y\|^2$, which satisfies the moment condition

$$\varlimsup_{N\to\infty} \frac{1}{\lambda_N^2} \log \int \exp(a\lambda_N^2 \varphi(Z_N, Y_N))\, d\pi^{\otimes N}$$

$$\leq \varlimsup_{N\to\infty} \frac{1}{\lambda_N^2} \log \int \exp(a\lambda_N^2 \langle\theta, Z_N\rangle)\, d\pi^{\otimes N}$$

$$= \varlimsup_{N\to\infty} \frac{N}{\lambda_N^2} \log \int \exp\left(\frac{\lambda_N}{\sqrt{N}}\langle a\theta, \widetilde{g}(x_1)\rangle\right) \pi(dx_1)$$

$$= \varlimsup_{N\to\infty} \frac{N}{\lambda_N^2} \left(1 + \frac{\lambda_N^2 a^2}{2N}\langle\theta, \Sigma_{11}\theta\rangle + o\left(\frac{\lambda_N^2}{N}\right)\right)$$

$$< \infty,$$

for any constant $a$. Since $\widetilde{C}_N$ is bounded in $N$, we obtain

$$\widetilde{\Lambda}(\theta) := \lim_{N\to\infty} \frac{1}{\lambda_N^2} \widetilde{\Lambda}_N(\theta)$$

$$= \lim_{N\to\infty} \frac{1}{\lambda_N^2} \log \int \exp \lambda_N^2 \varphi(Z_N, Y_N)\, d\pi^{\otimes N}$$

$$= \sup_{z,y}\{\varphi(z,y) - J(z,y)\}.$$

In order to maximize the right-hand side, write $(z,y)$ as $\Sigma(u,v)$, without loss of generality since $J$ is equal to $+\infty$ out of the range of $\Sigma$. Now

$$\sup_{z,y}\{\varphi(z,y) - J(z,y)\} = \sup_{u,v}\{\langle\theta, \Sigma_{11}u + \Sigma_{12}v\rangle - \tfrac{1}{2}\|\Sigma_{21}u + \Sigma_{22}v\|^2$$
$$- \tfrac{1}{2}(\langle u, \Sigma_{11}u\rangle + \langle v, \Sigma_{22}v\rangle + 2\langle u, \Sigma_{12}v\rangle)\}.$$

The function to be maximized is concave in $(u,v)$ and it is immediately checked that $(-\theta, (I+\Sigma_{22})^{-1}\Sigma_{21}\theta)$ is a stationary point. Substituting this back into the above expression, we finally arrive at

$$\widetilde{\Lambda}(\theta) = \tfrac{1}{2}\langle\theta, B\theta\rangle,$$

where $B = \Sigma_{11} - \Sigma_{12}(I+\Sigma_{22})^{-1}\Sigma_{21}$. In order to apply the Gartner–Ellis theorem, we need only to check that $B$ is nonnegative definite and apply

Lemma 4. Set $A = \Sigma_{12}\Sigma_{22}^-$. Since $\mathrm{Ker}[\Sigma_{22}] = \mathrm{Ker}[\Sigma_{12}]$, we have $\Sigma_{12} = A\Sigma_{22}$. As a consequence,

$$\mathrm{Var}_\pi[g(x_i) - AH(x_i)] = \Sigma_{11} - \Sigma_{12}\Sigma_{22}^-\Sigma_{21} \geq 0.$$

Now consider the difference

$$D = \Sigma_{12}\Sigma_{22}^-\Sigma_{21} - \Sigma_{12}(I+\Sigma_{22})^{-1}\Sigma_{21} = \Sigma_{12}(\Sigma_{22}^- - (I+\Sigma_{22})^{-1})\Sigma_{21},$$

and notice that the matrix between round brackets is nonnegative definite on $\mathrm{Ran}\,\Sigma_{22}$. But since $\mathrm{Ran}[\Sigma_{21}] \subset \mathrm{Ran}\,\Sigma_{22}$ (as a consequence of the inclusion $\mathrm{Ker}\,\Sigma_{22} \subset \mathrm{Ker}\,\Sigma_{12}$), $D$ is nonnegative definite, and, hence, so is $B$. The explicit estimate in Proposition 5 follows by taking $c = \inf\{\widetilde{\Lambda}^*(z) : z \notin B_1\} > 0$, where $B_1$ is the unit sphere in $\mathbb{R}^m$. $\square$

The results of this appendix can be directly applied to the sequence of densities $\pi_N$ defined in (4) by setting $m = 0$ and $\mu(dx) = \exp\{U(x)\}\,dx$.

## APPENDIX B.

Let $\mathcal{D}$ be the set of monomials in the derivatives of $H$ and $U$. By assumption (HP) functions in $\mathcal{D}$ are bounded. The following lemma is the result of a tedious but a straightforward computation, whose details are omitted.

LEMMA 6. For $h = 0, 1, 2, \ldots,$

$$\frac{d^h}{d\sigma^h}G_{\sigma,N}(x,W) = N\sum_{k=0}^{h+2}\sigma^k P_k(E_N\rho_\ell(x)\varphi_\ell(Y_\sigma)W^{r_\ell}; \ell = 1,\ldots,m_k), \tag{38}$$

for some integers $m_k$, where $P_k$ is a polynomial and $\rho_l, \varphi_l \in \mathcal{D}$. In particular, the derivatives $g_{k,N}(x,W) = \frac{1}{k!}\frac{d^k}{du^k}G_{u,N}(x,W)(0)$, for $k = 3,\ldots,6$, have the following explicit form:

$$g_{3,N}(x,W) = -\frac{N}{12}(E_N(3\psi_N''\psi_N'W + \psi_N'''W^3) \tag{39}$$
$$- 3E_N(H'W)E_N(H'\psi_N') - 3E_N(H''W^2)E_N(H'W)),$$

$$g_{4,N} = -\frac{N}{24}\Big(E_N(3\psi_N'''\psi_N'^2 + 3\psi_N''^2W^2 + 6\psi_N'''\psi_N'W^2 + \psi_N''''W^4)$$
$$- 3\{(E_N H'\psi_N')^2 + 2(E_N H''W^2)(E_N H'\psi_N') \tag{40}$$
$$+ (E_N H''W^2)^2\} + \delta_{4,N}\Big),$$

$$g_{5,N} = N\delta_{5,N}, \tag{41}$$

*and*

$$g_{6,N} = -\frac{N}{1440}$$
$$\times \Big\{ E_N(45\psi_N''^2\psi_N'^2 + 60\psi_N'''\psi_N'^3 + 90\psi_N'''\psi_N''\psi_N'W^2$$
$$+ 180\psi_N'''\psi_N''\psi_N' + 45(\psi_N''')^2 W^4$$
$$+ 180\psi_N''''\psi_N'^2 W^2 + 60\psi_N''''\psi_N'' W^4$$
$$+ 60\psi_N'''''\psi_N' W^4 + 4\psi_N''''''W^6)$$
$$- [90(E_N H'\psi_N'\psi_N'')(E_N H'\psi_N') + 180(E_N H''\psi_N'^2)(E_N H'\psi_N')$$
$$+ 90(E_N \psi_N''' H'W^2)(E_N H'\psi_N')$$
$$+ 90(E_N H'\psi_N''\psi_N')(E_N H''W^2)$$
$$+ 180(E_N H''\psi_N'')(E_N H'\psi_N') + 90(E_N \psi_N''' H'W^2)(E_N H''W^2)$$
$$+ 360(E_N H'''\psi_N'W^2)(E_N H'\psi_N')$$
$$+ 180(E_N H''W^2)(E_N H''\psi_N'^2)$$
$$+ 180(E_N H''W^2)(E_N H''\psi_N''W^2)$$
$$+ 60(E_N H''''W^4)(E_N H'\psi_N')$$
$$+ 360(E_N H''W^2)(E_N H'''\psi_N'W^2)$$
$$+ 60(E_N H''''W^4)(E_N H''W^2)]$$
$$+ [45(E_N H'^2)(E_N H'\psi_N')^2 + 90(E_N H'^2)(E_N H''W^2)(E_N H'\psi_N')$$
$$+ 45(E_N H'^2)(E_N H''W^2)^2]$$
$$+ \delta_{6,N} \Big\},$$

(42)

where $\delta_{4,N}, \delta_{5,N}$ and $\delta_{6,N}$ are sums of monomials in empirical averages of the type (15) and (18) and each of them has at least a factor with an odd value of $l$.

LEMMA 7. *Set*

$$\tau^2 = \tfrac{1}{144}\Big\{ 9\pi(\psi''^2(X)\psi'^2(X)) + 18\pi(\psi'(X)\psi''(X)\psi'''(X)) + 15\pi(\psi'''^2(X))$$
$$- 18\pi(H''(X) + H'(X)\psi'(X))\pi(H'(X)(\psi'''(X) + \psi'(X)\psi''(X)))$$
$$+ 9\pi(H'^2(X))(\pi(H''(X)) + \pi(H'(X)\psi'(X)))^2 \Big\}$$
$$=: F_3(\pi(\mathbf{r}_3(X)))$$

(43)

*for some polynomial $F_3$ and some vector $\mathbf{r}_3$ with components in $\mathcal{D}$.*

Then for any $N$ and $\varepsilon_N > 0$,

$$
\sup_u \left| \mathbb{P}\left( \frac{N^{-1/2} g_{3,N}(x, W)}{\tau} \le u \right) - \Phi(u) \right|
$$
(44)
$$
\le C\left( \frac{1}{\sqrt{N}} + \frac{1}{\varepsilon_N^2 N} \right) + h_\tau(F_3(E_N \mathbf{r}_3(x))) + \frac{\varepsilon_N}{\sqrt{2\pi}},
$$

where $h_\tau(x) = |1 \vee \frac{\sqrt{x}}{\tau}||1 - \frac{\tau}{\sqrt{x}}|$ is a continuous Lipschitz function vanishing at $\tau^2$.

PROOF. Let us define

$$
X_N = -\frac{\sqrt{N}}{12}\Big\{ 3E_N(\psi_N'' \psi_N'(x) W) + E_N(\psi_N'''(x) W^3) \\
- 3E_N(H' \psi_N'(x)) E_N(H'(x) W) - 3E_N H''(x) E_N(H'(x) W) \Big\}
$$

and

$$
Y_N = \frac{3\sqrt{N}}{12} E_N(H''(x)(W^2 - 1)) E_N(H'(x) W).
$$

From the expression of $g_{3,N}$ given in (39), we find that

$$
\frac{1}{\sqrt{N}} g_{3,N}(x, W) = X_N + Y_N.
$$

The term $Y_N$ has zero mean, and we bound its variance as follows:

(45)	$$\mathbb{E} Y_N{}^2 = \frac{9}{144 N^3} \sum_{i,j} H''^2(x_i) H'^2(x_j) \mathbb{E}((W_i^2 - 1)^2 W_j^2) \le \frac{C}{N}.$$

The expression $X_N$ is a sum of independent random variables, whose mean under the measure $\mathbb{P}$ is zero. We compute its variance $\tau_N^2$ directly as follows:

$$
\tau_N^2 = \frac{1}{144}\Big\{ E_N(9\psi_N''^2(x)\psi_N'^2(x) + 18\psi_N'(x)\psi_N''(x)\psi_N'''(x) + 15\psi_N'''^2(x)) \\
- 18 E_N[H'(x)\psi_N'(x) + H''(x)] \\
\times E_N[\psi_N'(x)\psi_N''(x)H'(x) + \psi_N'''(x)H'(x)] \\
+ 9 E_N H'^2(x)(E_N H''(x))^2 \\
+ 18 E_N H''(x) E_N(H'(x)\psi_N'(x)) E_N H'^2(x) \\
+ 9 E_N H'^2(x)(E_N(H'(x)\psi_N'(x)))^2
$$
(46)
$$
+ \frac{1}{N}\Big[ -36 E_N(\psi_N''(x)\psi_N'(x) H''(x) H'(x))
$$

$$- 48 E_N(\psi_N'''(x) H''(x) H'(x))$$
$$- 12 E_N(\psi_N'(x) \psi_N''(x) H'(x) H''(x))$$
$$- 48 E_N(\psi_N'''(x) H''(x) H'(x))$$
$$+ 36 E_N(H'(x) \psi_N'(x)) E_N(H''(x) H'^2(x))$$
$$+ 18 E_N(H'^2(x) H''(x)) E_N H''(x) + 36 E_N(H''^2(x) H'^2(x))\big]$$
$$+ 72 \frac{1}{N^2} E_N(H''^2(x) H'^2(x)) \Big\}.$$

By inserting into the above terms the explicit formula for $\psi_N$ given in (16), expanding the products and rearranging terms, we get the representation $\tau_N^2 = F_3(E_N \mathbf{r}_3(x))$. By replacing the vector of empirical averages $E_N \mathbf{r}_3(x)$ with that of expected values w.r.t. $\pi$, the expression (43) is obtained.

Next, setting $u = v \frac{\tau_N}{\tau}$, we obtain

$$\sup_u \left| \mathbb{P}\left(\frac{X_N}{\tau} \leq u\right) - \Phi(u) \right|$$
$$\leq \sup_v \left| \mathbb{P}\left(\frac{X_N}{\tau_N} \leq v\right) - \Phi(v) \right| + \sup_v \left| \Phi\left(v \frac{\tau_N}{\tau}\right) - \Phi(v) \right|$$
$$\leq \sup_v \left| \mathbb{P}\left(\frac{X_N}{\tau_N} \leq v\right) - \Phi(v) \right| + 1 \vee \left(\frac{\tau_N}{\tau}\right) \cdot \left| 1 - \left(\frac{\tau}{\tau_N}\right) \right|,$$

where the last line has been obtained by a straightforward Lipschitz estimate.

By using the formula given in [11], Lemma 1.9, page 20, again and the above estimate

$$\sup_u |\mathbb{P}(A_N \leq u) - \Phi(u)|$$
$$= \sup_u \left| \mathbb{P}\left(\frac{X_N + Y_N}{\tau} \leq u\right) - \Phi(u) \right|$$
$$\leq \sup_u \left| \mathbb{P}\left(\frac{X_N}{\tau_N} \leq u\right) - \Phi(u) \right| + \mathbb{P}(|Y_N| > \varepsilon_N \tau) + \frac{\varepsilon_N}{\sqrt{2\pi}}$$

and by means of Esseen's inequality ([11], Theorem 5.4, page 149) for $X_N/\tau_N$, Chebyshev's inequality and the estimate (45) for $Y_N$, we arrive at

$$\sup_u |\mathbb{P}(A_N \leq u) - \Phi(u)|$$
$$\leq \frac{1}{\sqrt{N}} \frac{C}{\tau_N^3} \Big\{ E_N |\psi_N''(x) \psi_N'(x)|^3 + E_N |\psi_N'''(x)|^3$$
$$+ E_N |H'(x)|^3 (E_N(|H'(x) \psi_N'(x)|^3 + |H''(x)|^3)) \Big\}$$

$$+ 1 \vee \left(\frac{\tau_N}{\tau}\right) \cdot \left|1 - \left(\frac{\tau}{\tau_N}\right)\right| + \frac{C}{N\tau^2 \varepsilon_N^2} E_N(H''^2(x)) E_N(H'^2(x)) + \frac{\varepsilon_N}{\sqrt{2\pi}},$$

from which the estimate (44) is obtained. □

REMARK 8. It is worth noting that when $H = 0$, that is, the target distribution has independent components, an easy integration by parts yields

$$\begin{aligned}\tau^2 &= \tfrac{1}{144}\{9\pi(\psi''^2(X)\psi'^2(X)) + 18\pi(\psi'(X)\psi''(X)\psi'''(X)) + 15\pi(\psi'''^2(X))\} \\ &= \tfrac{1}{48}\{5\pi(\psi'''^2(X)) - 3\pi(\psi''^3(X))\},\end{aligned}$$

which coincides with the constant $J^2$ appearing in the paper [14].

LEMMA 9. *Let $F:\mathbb{R}^m \to \mathbb{R}$ be a polynomial and $r_h:\mathbb{R}^2 \to \mathbb{R}$, $h = 1, \ldots, m$, be of the form $r_h(x_i, W_i) = b_h(x_i) W_i^{\beta_h}$, where $b_h$ belongs to $\mathcal{D}$. Define the vector $\mathbb{E}\mathbf{r}$ with the components in $\mathcal{D}$ by $(\mathbb{E}\mathbf{r})_h(x_i) = \mathbb{E}\{\mathbf{r}_h(x_i, W_i)\}$. Then for any $0 \leq \gamma < 1/2$ and $\varepsilon > 0$,*

$$\mathbb{P}[N^\gamma |F(E_N \mathbf{r}(x, W)) - F(\pi((\mathbb{E}\mathbf{r})(X)))| > \varepsilon] \leq \frac{C}{N^{(1-2\gamma)} \varepsilon^2}$$

*holds for all $x \in \widehat{F}_N(\varepsilon)$, where*

$$\widehat{F}_N(\varepsilon) = \{x : |E_N(\mathbb{E}\mathbf{r})(x) - \pi((\mathbb{E}\mathbf{r})(X))| < \varepsilon N^{-\gamma}/2K\},$$

*and $K$ is a local Lipschitz constant for $F$ in a neighborhood of the point $\pi((\mathbb{E}\mathbf{r})(X))$.*

PROOF. Let us notice that, when $x \in \widehat{F}_N(\varepsilon)$, we have

$$\begin{aligned}&\mathbb{P}(N^\gamma |F(E_N \mathbf{r}(x, W)) - F(\pi((\mathbb{E}\mathbf{r})(X)))| > \varepsilon) \\ &\quad \leq \mathbb{P}(N^\gamma |F(E_N(\mathbb{E}\mathbf{r})(x)) - F(E_N \mathbf{r}(x, W))| > \varepsilon/2).\end{aligned}$$

Let us consider a generic monomial appearing in $F(v_1, \ldots, v_m)$, which will be of the form $\prod_{h=1}^m v_h^{\alpha_h}$. By simple algebraic manipulations,

$$\prod_{h=1}^m v_h^{\alpha_h} - \prod_{h=1}^m u_h^{\alpha_h} = \sum_{(l_1, \ldots, l_m): l_1 + \cdots + l_m > 0} \prod_{h=1}^m \binom{\alpha_h}{l_h} (v_h - u_h)^{l_h} u_h^{\alpha_h - l_h}.$$

Now substitute the empirical average $E_N r_h(x, W)$ into $v_h$ and its centering $E_N(\mathbb{E}\mathbf{r})_h(x)$ into $u_h$. Denoting by $\mathbf{s} = \mathbf{r} - \mathbb{E}\mathbf{r}$, the above expression becomes

$$\sum_{(l_1, \ldots, l_m): l_1 + \cdots + l_m > 0} \prod_{h=1}^m \binom{\alpha_h}{l_h} (E_N s_h(x, W))^{l_h} E_N(\mathbb{E}\mathbf{r})_h(x)^{\alpha_h - l_h}.$$

OPTIMAL SCALING OF MALA 23We proceed to bound the second moment of each term of the above sum in the following way. The term $|\mathbb{E}\mathbf{r}|$ is bounded by a constant so we are left to bound the second moment

$$M_{h_1,\ldots,h_k}(x) = \mathbb{E}[(E_N s_{h_1}(x,W))^{\alpha_1} \cdots (E_N s_{h_k}(x,W))^{\alpha_k}]^2,$$

where $s_h(x_i, W_i) = b_h(x_i) Z_i^{(h)}$ with $Z_i^{(h)} = W_i^{\alpha_h} - \mathbb{E} W_i^{\alpha_h}$. By using next Lemma 10, we finally get the bound

$$M_{h_1,\ldots,h_k}(x) \leq \frac{C}{N}.$$

The proof is complete by an application of Chebyshev's inequality. □

LEMMA 10. *Let $(\mathbf{Z}_i : i = 1, \ldots, N)$ be i.i.d. centered $r$-dimensional random vectors. For any $j = 1, \ldots, r$ define $Y_i^{(j)} = b^{(j)}(x_i) Z_i^{(j)}$. Then for any $\alpha_j > 0$, $j = 1, \ldots, r$, such that $\sum_{j=1}^r \alpha_j = k$, it holds*

(47)
$$\mathbb{E}\left(\prod_{j=1}^r \left(\frac{1}{N} \sum_{i=1}^N Y_i^{(j)}\right)^{\alpha_j}\right)^2$$
$$\leq \frac{1}{N^k} \sum_{m=1}^k \frac{1}{N^{k-m}} \sum_{|\mathcal{P}|=m} \left(\frac{1}{N} \sum_{h_1=1}^N b^{A_1}(x_{h_1})\right) \cdots \left(\frac{1}{N} \sum_{h_m=1}^N b^{A_m}(x_{h_m})\right),$$

*where $b^{A_k}(x) = \mathbb{E} \prod_{j \in A_k} |b^{(j)}(x) Z_1^{(j)}|$ and the sum is taken over partitions $\mathcal{P} = \{A_1, \ldots, A_m\}$ of the set of repeated indices $I = \{1, \ldots, 1, 2, \ldots, 2, \ldots, r, \ldots, r\}$ (where "1" is repeated $2\alpha_1$ times, $\ldots$, "$r$" is repeated $2\alpha_r$ times) such that each $A_s$ contains at least two elements of $I$.*

PROOF. Begin by writing

$$\mathbb{E}\left[\prod_{j=1}^r \left(\frac{1}{N} \sum_{i=1}^N Y_i^{(j)}\right)^{\alpha_j}\right]^2$$

(48)
$$= \frac{1}{N^{2k}} \sum_{i_1=1}^N \cdots \sum_{i_k=1}^N \sum_{s_1=1}^N \cdots \sum_{s_k=1}^N \mathbb{E}(Y_{i_1}^{(1)} \cdots Y_{i_{\alpha_1}}^{(1)} \cdots Y_{i_{\alpha_1 + \cdots + \alpha_{r-1}}}^{(r)} \cdots Y_{i_k}^{(r)}$$
$$\times Y_{s_1}^{(1)} \cdots Y_{s_{\alpha_1}}^{(1)} \cdots Y_{s_{\alpha_1 + \cdots + \alpha_{r-1}}}^{(r)} \cdots Y_{s_k}^{(r)}).$$

A summand in the last expression is zero as soon as there exists an index $(i_1, \ldots, i_k, s_1, \ldots, s_k)$ whose value is *not* repeated by another. This follows by the independence and zero mean property of the $Y_i^{(j)}$. Another way of rearranging this sum is therefore as follows: partition the set $I$ of the upper indices of the formula (48) into a finite union $I = A_1 \cup \cdots \cup A_m$, where $|A_s| \geq 2$ for each $s$. We write $Y_i^{A_k} = \prod_{j \in A_k} Y_i^{(j)}$ to simplify notation. Then the sum on the left-hand side is bounded above in absolute value by

$$
(49) \qquad \sum_{m=1}^{k} \sum_{|\mathcal{P}|=m} \underbrace{\sum_{h_1=1}^{N} \cdots \sum_{h_m=1}^{N}}_{h_i \neq h_k \text{ if } k \neq i} \mathbb{E}|Y_{h_1}^{A_1}| \cdots |Y_{h_m}^{A_m}|.
$$

Since the sum is over nonrepeating indices $h_1, \ldots, h_m$, we have, by independence, $\mathbb{E}|Y_{h_1}^{A_1}| \cdots |Y_{h_m}^{A_m}| = b^{A_1}(x_{h_1}) \cdots b^{A_m}(x_{h_m})$. Now the summand in (49) is positive, so we can bound the sum from above by a sum over all (possibly repeating) indices $h_1, \ldots, h_m$, and after rearranging the sum, we obtain (47). □

LEMMA 11. *It holds that*

$$
(50) \qquad \mathbb{P}(|B_N| \geq \varepsilon_N) = \mathbb{P}\left(\frac{|g_{4,N}(x,W)|}{N^{2/3}} \geq \ell^{-1} \tau \varepsilon_N\right) \leq \frac{C}{N^{1/3} \varepsilon_N^2},
$$

$$
(51) \qquad \mathbb{P}(|C_N| \geq \varepsilon_N) = \mathbb{P}\left(\frac{|g_{5,N}(x,W)|}{N^{5/6}} \geq \ell^{-2} \tau \varepsilon_N\right) \leq \frac{C}{N^{2/3} \varepsilon_N^2},
$$

$$
(52) \qquad \mathbb{P}(|D_N| \geq \varepsilon_N) = \mathbb{P}\left(\left|\frac{g_{6,N}(x,W)}{N} + \frac{\tau^2}{2}\right| \geq \ell^{-3} \tau \varepsilon_N\right) \leq \frac{C}{N \varepsilon_N^2},
$$

*for $x \in \widehat{F}_{N,k}(\varepsilon_N)$, where*

$$
\widehat{F}_{N,k}(\varepsilon_N) = \left\{x : |E_N \mathbb{E}\mathbf{r}_k(x) - \pi(\mathbb{E}\mathbf{r}_k(X))| \leq \frac{\tau \varepsilon_N}{2K} \ell^{3-k} N^{k/6-1}\right\},
$$

*for $k = 4, 5, 6$, where $K$ is the smallest of the local Lipschitz constants for $F_k$ at $\pi(\mathbb{E}\mathbf{r}_k(X))$, for $k = 4, 5, 6$.*

PROOF. By Lemma 6 we have $g_{k,N}(x,W) = NF_k(E_N \mathbf{r}_k(x,W))$ for $k = 4, 5, 6$. The vectors $\mathbf{r}_k$ and polynomials $F_k$ are of the type required by Lemma 9. In order to compute $F_k(\pi(\mathbb{E}\mathbf{r}_k(X)))$ for $k = 4, 5, 6$ we need to replace in (40)–(42), of Lemma 6 the empirical averages with expectations with respect to $\pi \times \mathbb{P}$. By a straightforward computation,

$$
\begin{aligned}
F_4(\pi(\mathbb{E}\mathbf{r}_4(X))) &= -\tfrac{1}{24}\Big\{[3E(\psi''(X)\psi'^2(X)) + 3E(\psi'^2(X)) \\
&\qquad + 6E(\psi'''(X)\psi'(X)) + 3E(\psi''''(X))] \\
&\qquad - (E[H'(X)\psi'(X) + H''(X)])^2\Big\} \\
&= -\tfrac{1}{24}\left[3c \int_{-\infty}^{+\infty} (e^\psi \psi'')''(x)\,dx - \left(c \int_{-\infty}^{+\infty} (e^\psi H')'(x)\,dx\right)^2\right],
\end{aligned}
$$

since $X$ has the density $\pi(x) = ce^{\psi(x)}$ with $c = e^{-K(\theta_*)}$. Since by assumption (HP) both $(e^\psi \psi'')'(x)$ and $e^{\psi(x)} H'(x)$ are of the form $f(x)e^{\psi(x)}$ with $f$ bounded and $\psi(x) \to -\infty$ as $|x| \to +\infty$, the right-hand side of the previous expression is zero.

Next $F_5(\pi(\mathbb{E}\mathbf{r}_5(X))) = 0$, since each monomial in $\mathbf{r}_5$ contains at least one factor which is an odd power of $W$, hence, it has mean zero. Finally,

$$
\begin{aligned}
&F_6(\pi(\mathbb{E}\mathbf{r}_6(X))) \\
&= -\tfrac{1}{1440}\Big\{45E(\psi''^2(X)\psi'^2(X)) + 60E(\psi'''(X)\psi'^3(X)) \\
&\qquad + 270E(\psi'''(X)\psi''(X)\psi'(X)) + 135E(\psi'''^2(X)) \\
&\qquad + 180E(\psi''''(X)\psi'^2(X)) + 180E(\psi''''(X)\psi''(X)) \\
&\qquad + 180E(\psi'''''(X)\psi'(X)) + 60E(\psi''''''(X)) \\
&\qquad - 90[E(H''(X) + H'(X)\psi'(X)) \\
&\qquad\qquad \times E(H'(X)(\psi'''(X) + \psi'(X)\psi''(X)))] \\
&\qquad - 90[2E(H''(X)\psi'^2(X)) + 2E(H''(X)\psi''(X)) \\
&\qquad\qquad + 4E(H'''(X)\psi'(X)) + 2E(H''''(X))] \\
&\qquad\qquad \times E(H''(X) + H'(X)\psi'(X)) \\
&\qquad + 45E(H'^2(X))(E(H''(X)) + E(H'(X)\psi'(X)))^2\Big\} \\
&= -\tfrac{1}{1440}\Big\{(45E(\psi''^2(X)\psi'(X)^2) \\
&\qquad + 90E(\psi'(X)\psi''(X)\psi'''(X)) + 75E(\psi'''^2(X))) \\
&\qquad - 90E(H''(X) + H'(X)\psi'(X)) \\
&\qquad + E(H'(X)(\psi'''(X) + \psi'(X)\psi''(X))) \\
&\qquad + 45E(H'^2(X))(E(H''(X)) + E(H'(X)\psi'(X)))^2\Big\} \\
&\quad - \tfrac{60}{1440}\Big\{E(\psi'''(X)\psi'^3(X)) + 3E(\psi'(X)\psi''(X)\psi'''(X)) \\
&\qquad + E(\psi'''^2(X)) + 3E(\psi''''(X)\psi'^2(X)) \\
&\qquad + 3E(\psi''''(X)\psi''(X)) + 3E(\psi'''''(X)\psi'(X)) + E(\psi''''''(X))\Big\} \\
&\quad + \tfrac{180}{1440}\Big\{[E(H''(X)\psi'^2(X)) + E(H''(X)\psi''(X)) \\
&\qquad + 2E(H'''(X)\psi'(X)) + E(H''''(X))] \\
&\qquad \times [E(H''(X) + H'(X)\psi'(X))]\Big\},
\end{aligned}
$$

and this simplifies to

$$F_6(\mathbb{E}\mathbf{r}_6(X)) = -\frac{\tau^2}{2},$$

because the first term in curly braces equals $-\frac{\tau^2}{2}$ by (43), the second term is proportional to

$$E\Big(\psi'''(X)\psi'^3(X) + 3\psi'(X)\psi''(X)\psi'''(X) + \psi'''^2(X)3E(\psi''''(X)\psi''(X))$$
$$+ 3\psi''''(X)\psi'^2(X) + 3\psi'''''(X)\psi'(X) + \psi''''''(X)\Big)$$
$$= c\int_{-\infty}^{+\infty} (e^\psi \psi''')''' \, dx = 0,$$

and the third term in curly braces contains the multiplicative factor

$$E(H''(X) + H'(X)\psi'(X)) = c\int_{-\infty}^{+\infty}(e^\psi H')' \, dx = 0.$$

The last two displays equal zero by the same argument used before. Therefore,

$$\mathbb{P}\Big(\frac{|g_{4,N}(x,W)|}{N^{2/3}} \geq \ell^{-1}\tau\varepsilon_N\Big)$$
$$= \mathbb{P}(N^{1/3}|F_4(E_N\mathbf{r}_4(x,W))| \geq \ell^{-1}\tau\varepsilon_N),$$
$$\mathbb{P}\Big(\frac{|g_{5,N}|(x,W)}{N^{5/6}} \geq \ell^{-2}\tau\varepsilon_N\Big)$$
$$= \mathbb{P}(N^{1/6}|F_5(E_N\mathbf{r}_5(x,W))| \geq \ell^{-2}\tau\varepsilon_N),$$
$$\mathbb{P}\Big(\Big|\frac{g_{6,N}(x,W)}{N} + \frac{\tau^2}{2}\Big| \geq \ell^{-3}\tau\varepsilon_N\Big)$$
$$= \mathbb{P}(|F_6(E_N\mathbf{r}_6(x,W)) - F_6(\mathbb{E}\mathbf{r}_6(X))| \geq \ell^{-3}\tau\varepsilon_N)$$

so that the stated estimates follow directly from the previous lemma. □

LEMMA 12. *For $\sigma_N = \ell/N^{1/6}$, it holds*

$$(53) \quad \mathbb{P}\bigg[\bigg|\frac{1}{6!}\int_0^{\sigma_N}(\sigma_N - u)^6 \frac{d^7}{du^7}G_{u,N}(x,W)\,du\bigg| > \varepsilon_N\bigg] \leq \frac{C}{\varepsilon_N N^{1/6}}.$$

PROOF. By Markov's inequality and Lemma 6, we have

$$\mathbb{P}\bigg[\bigg|\frac{1}{6!}\int_0^{\sigma_N}(\sigma_N - u)^6 \frac{d^7}{du^7}G_{u,N}(x,W)\,du\bigg| > \varepsilon_N\bigg]$$
$$\leq \frac{1}{6!\varepsilon_N}\mathbb{E}\bigg|\int_0^{\sigma_N}(\sigma_N - u)^6 \frac{d^7}{du^7}G_{u,N}(x,W)\,du\bigg|$$

$$\leq \frac{1}{6!\varepsilon_N} \int_0^{\sigma_N} (\sigma_N - u)^6 \mathbb{E}\left|\frac{d^7}{du^7} G_{u,N}(x,W)\right| du$$

$$\leq \frac{1}{6!\varepsilon_N} \int_0^{\sigma_N} (\sigma_N - u)^6 N \mathbb{E}\left|\sum_{k=0}^{9} u^k P_k(E_N \rho_\ell(x) \varphi_\ell(Y_u) W^{r_\ell}; \ell = 1, \ldots, m)\right| du$$

$$\leq \frac{1}{6!\varepsilon_N} \int_0^{\sigma_N} (\sigma_N - u)^6 N \sum_{k=0}^{9} u^k \mathbb{E} |P_k(E_N \rho_\ell(x) \varphi_\ell(Y_u) W^{r_\ell}; \ell = 1, \ldots, m)| \, du,$$

$$\leq \frac{C}{\varepsilon_N} N {\sigma_N}^7 \leq \frac{C}{\varepsilon_N N^{1/6}}. \qquad \square$$

L. A. BREYER
UNIVERSITY OF LANCASTER
DEPARTMENT OF MATHEMATICS AND STATISTICS
FYLGE COLLEGE LANCASTER LA1 4YF
UNITED KINGDOM
E-MAIL: l.breyer@lancaster.ac.uk

M. PICCIONI
UNIVERSITY OF ROME LA SAPIENZA
DIPARTIMENTO DI MATEMATICA
PIAZZALE ALDO MORO 2
00185 ROME
ITALY
E-MAIL: piccioni@mat.uniroma1.it

S. SCARLATTI
UNIVERSITY "G. D'ANNUNZIO" CHIETI
DIPARTIMENTO DI SCIENZE
VIALE PINDARO 42
65127 PESCARA
ITALY
E-MAIL: scarlatt@sci.unich.it